\documentclass{amsart}
\usepackage{amssymb,amsxtra}

\usepackage{graphicx}
\usepackage{amscd}
\usepackage{amsmath}
\usepackage{amsfonts}
\usepackage{amssymb}
\title[iterating the Pimsner construction]{Iterating the Pimsner construction}
\author[deaconu]{Valentin Deaconu}
\address{Department of Mathematics and Statistics, University of Nevada, Reno NV
89557-0084, USA}
\email{vdeaconu@unr.edu}
\date{\today}

\begin{document}

\begin{abstract}
For $A$ a $C^*$-algebra, $E_1, E_2$ two  Hilbert bimodules over $A$, and a fixed isomorphism $\chi : E_1\otimes_AE_2\rightarrow E_2\otimes_AE_1$, we consider the problem of computing the $K$-theory of the Cuntz-Pimsner algebra ${\mathcal O}_{E_2\otimes_A{\mathcal O}_{E_1}}$ obtained by extending the scalars and by iterating the Pimsner construction.

The motivating examples are a commutative diagram of Douglas and Howe for the Toeplitz operators on the quarter plane, and the Toeplitz extensions associated by Pimsner and Voiculescu to compute the $K$-theory of a crossed product. The applications are for Hilbert bimodules arising from  rank two graphs and from commuting endomorphisms of abelian $C^*$-algebras.
\end{abstract}

\maketitle

\centerline{\S 0. {\small INTRODUCTION}}

\bigskip

In his seminal paper \cite{Pi},  Pimsner introduced a large class of $C^*$-algebras generalizing both the Cuntz-Krieger algebras and the crossed products by an automorphism. The central notion is that of a Hilbert bimodule or $C^*$-correspondence, which appeared also in the theory of subfactors. His construction was modified by Katsura (see \cite{Ka1}) to include bimodules defined from graphs with sinks, or more generally, from topological graphs.

In the same way that a crossed product by ${\mathbb Z}^2$ could be thought as an iterated crossed product by ${\mathbb Z}$, we iterate the Pimsner construction for two
``commuting'' Hilbert bimodules over a $C^*$-algebra.  A basic example comes from a rank two graph of Kumjian and Pask (see \cite{KP}), where the two bimodules correspond to the horizontal and vertical edges, and the commutation relation is given by the unique factorization property.  We get a particular case of the product systems defined by Fowler (see \cite{F2}), but our approach has the advantage that it generates some exact sequences of $K$-theory.

After some preliminaries about the  algebras ${\mathcal T}_E$ and ${\mathcal O}_E$, we describe how from a Hilbert module $E$ over $A$ and a map $A\to B$ we can get a Hilbert module over $B$, using a tensor product.  Given two Hilbert bimodules $E_1, E_2$ over $A$ and an isomorphism $\chi: E_1\otimes_AE_2\to E_2\otimes_AE_1$ we consider $E_2\otimes_A{\mathcal T}_{E_1}$ and $E_2\otimes_A{\mathcal O}_{E_1}$ as Hilbert bimodules over ${\mathcal T}_{E_1}$ and ${\mathcal O}_{E_1}$ and repeat the Pimsner construction. The last section deals with a $3\times 3$  diagram involving the
iterated Cuntz-Pimsner algebras, inspired from the work of Douglas and Howe (see \cite{DH}) and Pimsner and Voiculescu (see \cite{PV}). As a corollary, we obtain some exact sequences of $K$-theory, including information about the $K$-groups of the $C^*$-algebra ${\mathcal O}_{E_2\otimes_A{\mathcal O}_{E_1}}$. Several examples are considered, including commuting endomorphisms of abelian $C^*$-algebras. 

Since this paper was completed, we discovered that  J. Lindiarni and I. Raeburn (\cite{LR}) already used  the diagram which appears in our Lemma 4.2.

\newpage

\centerline{\S 1. {\small PRELIMINARIES}}

\bigskip

Recall that a (right) Hilbert $A$-module is a Banach space $E$ with a right action of a $C^*$-algebra $A$ and an $A$-valued inner product $\langle\cdot,\cdot\rangle:E\times E\rightarrow A$ linear in the second variable such that
\[\langle\xi,\eta a\rangle=\langle\xi,\eta\rangle a, \; \langle\xi,\eta\rangle=\langle\eta,\xi\rangle^*,\; \langle\xi,\xi\rangle\geq 0,\;||\xi||=||\langle\xi,\xi\rangle||^{1/2}.\]
A Hilbert module is called full if the closed linear span of the inner products coincides with $A$. We denote by ${\mathcal L}(E)$ the $C^*$-algebra of adjointable operators on $E$, and by $\theta_{\xi,\eta}\in{\mathcal L}(E)$ the rank one operator
\[\theta_{\xi,\eta}(\zeta)=\xi\langle\eta,\zeta\rangle.\]
The closed linear span of rank one operators is the ideal ${\mathcal K}(E)$. We have ${\mathcal L}(E)\cong M({\mathcal K}(E))$, the multiplier algebra. Also, ${\mathcal K}(E)$ can be identified with the balanced tensor product $E\otimes_AE^*$, where $E^*$ is the dual of $E$, a left Hilbert $A$-module.

A Hilbert bimodule over $A$ (sometimes called a $C^*$-correspondence from $A$ to $A$) is a Hilbert $A$-module with a left action of $A$ given by a homomorphism $\varphi:A\rightarrow {\mathcal L}(E)$. A Hilbert bimodule  is called  faithful if $\varphi$ is injective. The left action is nondegenerate if $\overline{\varphi(A)E}=E$.  For $n\ge 0$ we denote by $E^{\otimes n}$ the Hilbert bimodule obtained by taking the tensor product of $n$ copies of $E$, balanced over $A$ (for $n=0, E^{\otimes 0}=A$). Recall that for $n=2$, the inner product is given by \[\langle\xi\otimes\eta,\xi'\otimes\eta'\rangle=\langle\eta,\varphi(\langle\xi,\xi'\rangle)\eta'\rangle,\] and it is inductively defined for general $n$.

\medskip

{\bf 1.1 Definition}. A {\em Toeplitz representation} of a Hilbert bimodule $E$ over $A$ in a $C^*$-algebra $C$ is a pair $(\tau,\pi)$ with $\tau:E\rightarrow C$ a linear map and $\pi:A\rightarrow C$ a *-homomorphism, such that 
\[\tau(\xi a)=\tau(\xi)\pi(a), \;\tau(\xi)^*\tau(\eta)=\pi(\langle \xi,\eta\rangle), \;\tau(\varphi(a)\xi)=\pi(a)\tau(\xi).\]
Note that the first property actually follows from the second. Indeed,
\[||\tau(\xi a)-\tau(\xi)\pi(a)||^2=(\tau(\xi a)-\tau(\xi)\pi(a))^*(\tau(\xi a)-\tau(\xi)\pi(a))=\]\[=
\pi(\langle\xi a,\xi a\rangle )-\pi(a)^*\pi(\langle\xi, \xi a\rangle)-\pi(\langle\xi a, \xi\rangle)\pi(a)+\pi(a)^*\pi(\langle\xi,\xi\rangle)\pi(a)=0.\]The corresponding universal $C^*$-algebra is called the Toeplitz algebra of $E$, denoted by ${\mathcal T}_E$. If $E$ is full, then ${\mathcal T}_E$ is generated by elements $\tau^n(\xi)\tau^m(\eta)^*, m,n\ge 0$, where $\tau^0=\pi$  and for $n\geq 1, \tau^n(\xi_1\otimes ...\otimes\xi_n)=\tau(\xi_1)...\tau(\xi_n)$ is the extension of $\tau$ to $E^{\otimes n}$. If $E$ is also faithful, then $A\subset {\mathcal T}_E$.

There is a homomorphism $\psi:{\mathcal K}(E)\rightarrow C$ such that $\psi(\theta_{\xi,\eta})=\tau(\xi)\tau(\eta)^*.$
A representation $(\tau,\pi)$ is {\em Cuntz-Pimsner covariant} if $\pi(a)=\psi(\varphi(a))$ for all $a$ in the ideal \[I_E=\varphi^{-1}({\mathcal K}(E))\cap (\ker \varphi)^\perp.\] 
The Cuntz-Pimsner algebra ${\mathcal O}_E$ is universal with respect to the covariant representations, and it is a quotient of ${\mathcal T}_E$.

There is a gauge action of ${\mathbb T}$ on ${\mathcal T}_E$ and ${\mathcal O}_E$ defined by
\[z\cdot(\tau^n(\xi)\tau^m(\eta)^*)=z^{n-m}\tau^n(\xi)\tau^m(\eta)^*, z\in{\mathbb T},\]
and using the universal properties. The {\em core} ${\mathcal F}_E$ is the fixed point algebra ${\mathcal O}_E^{\mathbb T}$, generated by the union of the algebras ${\mathcal K}(E^{\otimes n})$.

The Toeplitz algebra ${\mathcal T}_E$ can be represented by creation operators $T_\xi(\eta)=\xi\otimes\eta$ on the Fock bimodule \[\ell^2(E)=\bigoplus_{n\geq 0}E^{\otimes n},\]  and there is an ideal $I={\mathcal K}(\ell^2(E)I_E)$ in ${\mathcal L}(\ell^2(E))$ such that 
\[0\rightarrow I\rightarrow {\mathcal T}_E\rightarrow {\mathcal O}_E\rightarrow 0\]is exact.
In particular, if ${\mathcal L}(E)={\mathcal K}(E)$ or $\varphi(A)\subset{\mathcal K}(E)$, then 
\[0\rightarrow {\mathcal K}(\ell^2(E))\rightarrow{\mathcal T}_E\rightarrow {\mathcal O}_E\rightarrow 0.\]
For more details about the algebras ${\mathcal F}_E, {\mathcal T}_E$, and ${\mathcal O}_E$ we refer to the original paper of Pimsner (\cite{Pi}) and to \cite{Ka1}.
 
If $E$ is finitely generated, then it has a basis $\{u_i\}$, in the sense that for all $\xi\in E$,
\[\xi=\sum_{i=1}^nu_i\langle u_i,\xi\rangle.\]
In this case, ${\mathcal L}(E)={\mathcal K}(E)$, and  the Cuntz-Pimsner algebra ${\mathcal O}_E$ is generated by $S_i=S_{u_i}$ with relations
\[\sum_{i=1}^nS_iS_i^*=1,\; S_i^*S_j=\langle u_i,u_j\rangle,\; a\cdot S_j=\sum_{i=1}^nS_i\langle u_i,a\cdot u_j\rangle.\]
The Toeplitz algebra ${\mathcal T}_E$ is generated by $T_i=T_{u_i}$ which satisfy only the last two relations (see \cite{KPW}).

\medskip

{\bf 1.2 Examples}. 1) For $A={\mathbb C}$ and  $E=H$ a Hilbert space with orthonormal basis $\{\xi_i\}_{i\in I}$, the Toeplitz algebra ${\mathcal T}_H$ is generated by $\{S_i\}_{i\in I}$ satisfying $S_i^*S_j=\delta_{ij}\cdot 1,\; i,j\in I$. In ${\mathcal O}_H$ we also have $\displaystyle\sum_{i\in I} S_iS_i^*=1$ if the dimension of $H$ is finite.  In particular, for $E={\mathbb C}=A$ we get ${\mathcal T}_E\cong{\mathcal T}$, the classical Toeplitz algebra generated by the unilateral shift,  ${\mathcal O}_E\cong C({\mathbb T})$, the continuous functions on the unit circle, and ${\mathcal F}_E\cong {\mathbb C}$. For $E={\mathbb C}^n$, we get ${\mathcal T}_E\cong{\mathcal E}_n$, the Cuntz-Toeplitz algebra,  ${\mathcal O}_E\cong{\mathcal O}_n$, the Cuntz algebra, and  ${\mathcal F}_E\cong UHF(n^{\infty})$. For $H$ infinite dimensional and separable, ${\mathcal T}_H\cong{\mathcal O}_H\cong{\mathcal O}_{\infty}$.

2) Let $E=A^n$ with the usual structure:
\[\langle(a_1,...,a_n),(b_1,...,b_n)\rangle=\sum_{i=1}^n a_i^*b_i,\; (a_1,...,a_n)\cdot a=(a_1a,...,a_na),\;a\cdot(a_1,...,a_n)=(aa_1,...,aa_n).\]
We get ${\mathcal T}_E\cong A\otimes{\mathcal E}_n,\;{\mathcal O}_E\cong A\otimes {\mathcal O}_n, \;{\mathcal F}_E\cong A\otimes UHF(n^{\infty})$.

3) Let $\alpha :A\rightarrow A$ be an automorphism of a unital $C^*$-algebra, and let $E=A(\alpha)$ be the Hilbert bimodule obtained from $A$ with 
the usual inner product and right multiplication, and with left action $\varphi(a)x=\alpha(a)x$. Then ${\mathcal T}_E$ is isomorphic to the Toeplitz extension ${\mathcal T}_{\alpha}$ used by Pimsner and Voiculescu in {\cite {PV}}, and ${\mathcal O}_E$ is isomorphic to the crossed product $A\rtimes_{\alpha}{\mathbb Z}$.
Indeed, let $\hat{a}$ denote the element in $E$ obtained from $a\in A$. Then $S=\tau(\hat{1})$ is an isometry in any unital Toeplitz representation $(\tau,\pi)$, since $\langle\hat{1},\hat{1}\rangle=1$, and it is an unitary in any unital covariant representation, since the rank one operator
$\theta_{\hat{1},\hat{1}}$ is the identity. We also have \[\pi(\alpha(a))=\pi(\langle\hat{1},\widehat{\alpha(a)}\rangle)=\tau(\hat{1})^*\tau(\widehat{\alpha(a)})=\tau(\hat{1})^*\tau(\varphi(a)\hat{1})=\tau(\hat{1})^*\pi(a)\tau(\hat{1})=S^*\pi(a)S,\] therefore ${\mathcal O}_E\cong A\rtimes_{\alpha}{\mathbb Z}$. In the paper mentioned above, the Toeplitz extension ${\mathcal T}_{\alpha}$ was defined as the $C^*$-subalgebra of $(A\rtimes_{\alpha}{\mathbb Z})\otimes{\mathcal  T}$ generated by $A\otimes 1$ and $u\otimes S_+$, where $\alpha(a)=uau^*$ and $S_+$ is the unilateral shift. It is easy to see that ${\mathcal  T}_{\alpha}\cong {\mathcal T}_E$ by the map which takes $a\otimes 1$ into $a$ and $u\otimes S_+$ into $S$. Since ${\mathcal K}(\ell^2(E))\cong A\otimes{\mathcal K}$,  we recover the short exact sequence
\[0\rightarrow A\otimes{\mathcal K}\rightarrow {\mathcal T}_{\alpha}\rightarrow A\rtimes_{\alpha}{\mathbb Z}\rightarrow 0.\]

4) Graph $C^*$-algebras. For an oriented countable graph $G=(G^0,G^1,r,s), \; C^*(G)$ is defined as the universal $C^*$-algebra generated by mutually orthogonal projections $\{p_v\}_{v\in G^0}$ and partial isometries $\{s_e\}_{e\in G^1}$ with orthogonal ranges such that $s_e^*s_e=p_{r(e)}, \;s_es_e^*\leq p_{s(e)}$ and \[(*)\;\;p_v=\sum_{s(e)=v}s_es_e^*\;\mbox{if}\;\;0<|s^{-1}(v)|<\infty.\]
We can set $A=C_0(G^0)$, and denote by $E$ the Hilbert module we obtain after we complete $C_c(G^1)$ in the norm given by the inner product
\[\langle\xi,\eta\rangle(v)=\sum_{r(e)=v}\overline{\xi(e)}\eta(e)\]
with the right action  defined by $(\xi f)(e)=\xi(e)f(r(e))$. The left action is defined by \[\varphi:A\rightarrow {\mathcal L}(E),\;  \varphi(f)\xi(e)=f(s(e))\xi(e),\; \xi\in C_c(G^1)\subset E.\] We have
\[\varphi^{-1}({\mathcal K}(E))=C_0(\{v\in G^0: |s^{-1}(v)|<\infty\}),\;\ker(\varphi)=C_0(\{v\in G^0: |s^{-1}(v)|=0\}),\]
hence $I_E=C_0(\{v\in G^0:0<|s^{-1}(v)|<\infty\})$. Define 
\[\pi(f)=\sum_{v\in G^0}f(v)p_v, \;f\in C_0(G^0), \;\;\tau(\xi)=\sum_{e\in G^1}\xi(e)s_e, \;\xi\in C_c(G^1)\subset E.\] The pair $(\tau, \pi)$ is a Toeplitz representation into $C^*(G)$ iff $s_e^*s_e=p_{r(e)}$ and  $s_es_e^*\leq p_{s(e)}$, which is covariant iff $(*)$ is satisfied. This proves that $C^*(G)\cong {\mathcal O}_E$ (see \cite{Ka1}). The core ${\mathcal F}_E$ is an AF-algebra. For an irreducible oriented finite graph with no sinks, we obtain the Cuntz-Krieger algebras ${\mathcal O}_A$ as   Cuntz-Pimsner algebras.

5) For a $C^*$-algebra $A$ and an injective unital endomorphism $\alpha\in End(A)$ such that there is a conditional expectation $P$ onto the range $\alpha(A)$, one can define a Hilbert bimodule $E=A(\alpha,P)$, using the transfer operator $L=\alpha^{-1}\circ P:A\to A$ as in \cite{EV}. We complete the vector space $A$ with respect to the inner product $\langle\xi,\eta\rangle=L(\xi^*\eta)$, and define the right and left multiplications by $\xi\cdot a=\xi\alpha(a),\; a\cdot\xi=a\xi$. We have 
\[\langle\xi,\eta\cdot a\rangle=\langle\xi,\eta\alpha(a)\rangle=L(\xi^*\eta\alpha(a))=\alpha^{-1}(P(\xi^*\eta\alpha(a)))=\alpha^{-1}(P(\xi^*\eta)\alpha(a))=\alpha^{-1}(P(\xi^*\eta))a=\langle\xi,\eta\rangle a.\]
For $\alpha$ an automorphism and $P=id$, the corresponding $C^*$-algebra ${\mathcal O}_E$ is isomorphic to the crossed product $A\rtimes_{\alpha^{-1}}{\mathbb Z}$. Indeed, let $F$ be the Hilbert bimodule $A(\alpha^{-1})$ with the structure as in example 3. The map $\alpha^{-1}:A\to A$ induces an isomorphism of Hilbert bimodules $h:E\to F$. For $A=C(X)$ with $X$ compact, and $\alpha$ induced by a surjective local homeomorphism $\sigma:X\rightarrow X$, we take
\[(Pf)(x)=\frac{1}{\nu(x)}\sum_{\sigma(y)=\sigma(x)}f(y),\]
where $\nu(x)$ is the number of elements in the fiber $\sigma^{-1}(x)$. It was proved in \cite{D} that the corresponding algebra ${\mathcal O}_{A(\alpha,P)}$ is isomorphic to $C^*(\Gamma(\sigma))$, where $\Gamma(\sigma)$ is the Renault groupoid
\[\Gamma(\sigma)=\{(x,p-q,y)\in X\times{\mathbb Z}\times X \;\mid\; \sigma^p(x)=\sigma^q(y)\}.\]

\bigskip

\centerline{\S 2. {\small EXTENDING THE SCALARS}}

\bigskip

Let $E$ be a Hilbert module over $A$ and let $\rho:A\to B$ be  a $C^*$-algebra homomorphism (we will be interested mostly in the case when $\rho$ is an inclusion). Then $B$ is a left $A$-module with multiplication $a\cdot b=\rho(a)b$, and $E\otimes_AB$ becomes a Hilbert module over $B$, with the inner product  given by
\[\langle\xi_1\otimes b_1,\xi_2\otimes b_2\rangle=b_1^*\rho(\langle\xi_1,\xi_2\rangle) b_2\]
and right multiplication $(\xi\otimes b_1)\cdot b_2=\xi\otimes b_1b_2$. We have ${\mathcal K}(E\otimes_AB)\cong (E\otimes_AB)\otimes_B(E\otimes_AB)^*\cong E\otimes_A B\otimes_A E^*$, which is strongly Morita equivalent to $B$ in the case $E$ is full. Also, we have an inclusion ${\mathcal K}(E)\subset{\mathcal K}(E\otimes_AB)$ for $B$ unital, given by $\xi\otimes\eta^*\mapsto\xi\otimes1\otimes\eta^*$. If $E$ is a Hilbert bimodule over $A$, and if there is a *-morphism
$B\rightarrow {\mathcal L}(E)$ which extends the left multiplication of $A$ on $E$, then $E\otimes_AB$ becomes a Hilbert bimodule over $B$, and one can form the tensor powers
$(E\otimes_AB)^{\otimes n}$. Assuming that the left action of $B$ on $E\otimes_AB$ is nondegenerate, we get
\[(E\otimes_AB)^{\otimes n}\cong E^{\otimes n}\otimes_AB.\]
In particular,   the Toeplitz algebra
${\mathcal T}_{E\otimes_AB}$ is represented on $\ell^2(E\otimes_AB)\cong \ell^2(E)\otimes_AB$, and depends on the left multiplication of $B$. An interesting question is to relate the $C^*$-algebras ${\mathcal T}_{E\otimes_AB}$ and  ${\mathcal O}_{E\otimes_AB}$ to ${\mathcal T}_E$, ${\mathcal O}_E$, and $B$.

\medskip

{\bf 2.1 Example}. Let $A={\mathbb C}$, let $E=H$  be a separable infinite dimensional Hilbert space, and let $B$ be a separable unital $C^*$-algebra. Then  $H\otimes B$ is a Hilbert module over $B$, and  ${\mathcal K}(H\otimes B)\cong {\mathcal K}(H)\otimes B$. If $B$ is faithfully represented on $H$, then $H\otimes B$ becomes a Hilbert bimodule over $B$. Assuming, in addition, that the intersection of $B$ with  ${\mathcal K}(H)$ is trivial,  Kumjian (see \cite {K}) showed that ${\mathcal T}_{H\otimes B}\cong {\mathcal O}_{H\otimes B}$ is simple and purely infinite, with the same $K$-theory as $B$.

\medskip

{\bf 2.2 Example}. Let $A=C_0(X)$ and let $E$ be a Hilbert module  given by a continuous field of elementary $C^*$-algebras over $X$. Then for an abelian $C^*$-algebra $B$ containing $C_0(X)$, the tensor product $E\otimes_AB$ is obtained by a pull-back. In particular,
let $(G^0,G^1,r,s)$ be a (topological) graph, and consider a covering map $p:\tilde{G^0}\rightarrow G^0$ which gives an inclusion $A=C_0(G^0)\subset C_0(\tilde{G^0})=B$.
The Hilbert module $E\otimes_AB$, where $E$ is obtained from $C_c(G^1)$ as in example 4 \S 1,
is associated to a``pull-back  graph'' in which the set of vertices is $\tilde{G^0}$ and the set of edges is $\tilde{G^1}:=\{(x,e,y)\in \tilde{G^0}\times G^1\times \tilde{G^0}\mid s(e)=p(x), r(e)=p(y)\}$. The new range and source maps are $ s(x,e,y)=x, r(x,e,y)=y$.
It is known that ${\mathcal K}(E)$ is isomorphic to $C^*(R)$, where $R$ is the equivalence relation \[R=\{(e_1,e_2)\in G^1\times G^1\mid r(e_1)=r(e_2)\}.\] 
Then ${\mathcal K}(E\otimes_AB)$ is isomorphic to $C^*(p^*(R))$, where
\[p^*(R)=\{((x_1,e_1,y_1),(x_2,e_2,y_2))\in \tilde{G^1}\times\tilde{G^1}\mid p(y_1)=p(y_2)\}.\]

\medskip

{\bf 2.3 Example}. Let $\alpha:A\rightarrow A$ be an automorphism of a $C^*$-algebra $A$, which extends to $\tilde{\alpha}:B\rightarrow B$, where $A\subset B$. Consider $E=A(\alpha)$ as in example 1.2.3. Then $E\otimes_AB\cong B(\tilde{\alpha})$, ${\mathcal T}_{E\otimes_AB}\cong{\mathcal T}_{\tilde{\alpha}}$ and ${\mathcal O}_{E\otimes_AB}\cong B\rtimes_{\tilde{\alpha}}{\mathbb Z}$, which  contains ${\mathcal O}_E\cong A\rtimes_\alpha{\mathbb Z}$.

\medskip

{\bf 2.4 Example}. For a Hilbert bimodule $E$ over $A$, Pimsner used the Hilbert module
$E_{\infty}=E\otimes_A{\mathcal F}_E$ (see \cite{Pi}, section 2) in order to get an inclusion ${\mathcal T}_E\subset{\mathcal T}_{E_{\infty}}$, an isomorphism ${\mathcal O}_E\cong {\mathcal O}_{E_{\infty}}$, and a completely positive map $\phi:{\mathcal O}_E\to{\mathcal T}_{E_{\infty}}$ which is a cross-section to the quotient map ${\mathcal T}_{E_{\infty}}\to{\mathcal O}_E$.

\bigskip

\centerline{\S 3. {\small ITERATING THE PIMSNER CONSTRUCTION}}

\bigskip

Consider now two full finitely generated Hilbert $A$-bimodules $E_1$ and $E_2$ such that $A$ is unital and the left actions $\varphi_i:A\rightarrow {\mathcal L}(E_i)$ are injective and nondegenerate. We assume that  there is an isomorphism of Hilbert $A$-bimodules $\chi:E_1\otimes_AE_2\rightarrow E_2\otimes_AE_1$. This isomorphism should be understood as a kind of commutation relation. The most interesting cases are when $E_1$ and $E_2$ are independent, in the sense that no tensor power of one is isomorphic to the other.

Note that the isomorphism $\chi$ induces and isomorphism $E_1^*\otimes_AE_2^*\rightarrow E_2^*\otimes_AE_1^*$ because $(E_1\otimes_AE_2)^*\cong E_2^*\otimes_AE_1^*$.

Since $A\subset {\mathcal T}_{E_1}$, the tensor product $E_2\otimes_A{\mathcal T}_{E_1}$ becomes a Hilbert module over ${\mathcal T}_{E_1}$ as in \S 2, with the inner product  given by
$\langle\xi\otimes x,\eta\otimes y\rangle=x^*\langle\xi,\eta\rangle y,$
and the right multiplication by $(\xi\otimes x)y=\xi\otimes xy$. Since ${\mathcal T}_{E_1}$ is generated by $E_1$, in order to define a left multiplication of ${\mathcal T}_{E_1}$ on $E_2\otimes_A{\mathcal T}_{E_1}$, it is sufficient to define the left multiplication  by elements in $E_1$. This is done via the composition
\[E_1\otimes_AE_2\otimes_A{\mathcal T}_{E_1}\rightarrow E_2\otimes_AE_1\otimes_A{\mathcal T}_{E_1}\rightarrow E_2\otimes_A{\mathcal T}_{E_1},\]
where the first map is $\chi\otimes id$, and the second is given by absorbing $E_1$ into ${\mathcal T}_{E_1}$. For the adjoint, note that there is  a map
\[E_1^*\otimes_AE_2\otimes_AE_1\to E_1^*\otimes_AE_1\otimes_AE_2\to E_2,\]
where the first map is $id\otimes\chi^{-1}$, and the second is given by left multiplication with the inner product in $E_1$. We will denote by $E_2\otimes^{\chi}_A{\mathcal T}_{E_1}$ the bimodule obtained using this left multiplication, when $\chi$ is not understood. In the same way, we may consider the bimodule $E_1\otimes_A{\mathcal T}_{E_2}$ with the left multiplication induced by $\chi^{-1}$.

\medskip

{\bf 3.1 Lemma}. With the above structure, $E_2\otimes_A{\mathcal T}_{E_1}$ is a Hilbert bimodule over ${\mathcal T}_{E_1}$, and we can consider the Toeplitz algebra ${\mathcal T}_{E_2\otimes_A{\mathcal T}_{E_1}}$. We have 
\[{\mathcal T}_{E_2\otimes_A{\mathcal T}_{E_1}}\cong {\mathcal T}_{E_1\otimes_A{\mathcal T}_{E_2}}.\]

{\em Proof.}  Both algebras ${\mathcal T}_{E_2\otimes_A{\mathcal T}_{E_1}},  {\mathcal T}_{E_1\otimes_A{\mathcal T}_{E_2}}$ are represented on the Fock space $\ell^2(E_2)\otimes_A\ell^2(E_1)\cong \ell^2(E_1)\otimes_A\ell^2(E_2)$, and are generated by ${\mathcal T}_{E_1}$ and ${\mathcal T}_{E_2}$ with the commutation relation given by the isomorphism $\chi$. $\Box$

Similarly, we can construct the Hilbert bimodules $E_2\otimes_A{\mathcal O}_{E_1}$ and
$E_1\otimes_A{\mathcal O}_{E_2}$. We get

\medskip

{\bf 3.2 Lemma}.  With the above notation,
\[{\mathcal T}_{E_2\otimes_A{\mathcal O}_{E_1}}\cong {\mathcal O}_{E_1\otimes_A{\mathcal T}_{E_2}},
\;\;{\mathcal O}_{E_2\otimes_A{\mathcal O}_{E_1}}\cong {\mathcal O}_{E_1\otimes_A{\mathcal O}_{E_2}}.\]

{\em Proof.} The first two algebras are quotients of ${\mathcal T}_{E_2\otimes_A{\mathcal T}_{E_1}}\cong {\mathcal T}_{E_1\otimes_A{\mathcal T}_{E_2}}$ by the ideal generated by ${\mathcal K}(\ell^2(E_1))$, and the last two are quotients by the ideal generated by ${\mathcal K}(\ell^2(E_1))$ and ${\mathcal K}(\ell^2(E_2))$. $\Box$

Note that there is a gauge action of ${\mathbb T}^2$ on ${\mathcal O}_{E_2\otimes_A{\mathcal O}_{E_1}}\simeq {\mathcal O}_{E_1\otimes_A{\mathcal O}_{E_2}}$.

\medskip

{\bf 3.3 Remark}.  Given $E_1, E_2$ as above, we can define a product system $E=E^{\chi}$ of Hilbert bimodules over the semigroup $({\mathbb N}^2,+)$ (see \cite{F2}), as follows. Define the fibers by $E_{(m,n)}=E_1^{\otimes m}\otimes_A E_2^{\otimes n}$ for $(m,n)\in {\mathbb N}^2$, and the multiplication induced by the isomorphism $\chi$. It is easy to see that we get associativity, and therefore we may consider the $C^*$-algebras ${\mathcal T}_E, \;{\mathcal O}_E$ as defined by Fowler. Note that if we change $\chi$, the product system also changes. In particular, a path of isomorphisms $\chi_t$ will determine a family of product systems $E^t=E^{\chi_t}$.

Recall that a Toeplitz representation of a product system
$E$ in a $C^*$-algebra is obtained from a family of Toeplitz representations of the fibers, compatible with the product. In this particular case, it is sufficient to consider two Toeplitz representations
$(\tau_1,\pi), (\tau_2, \pi)$ of the generators $E_1$ and $E_2$, respectively, and to define $\tau(\xi_1\otimes ...\otimes\xi_m\otimes\eta_1\otimes ...\otimes\eta_n)=\tau_1(\xi_1)...\tau_1(\xi_m)\tau_2(\eta_1)...\tau_2(\eta_n)$ for $\xi_i\in E_1, i=1,...,m$ and  $\eta_j\in E_2, j=1,...,n$.
The Toeplitz algebra ${\mathcal T}_E$ is represented on the Fock space $\ell^2(E)\cong
\ell^2(E_1)\otimes_A\ell^2(E_2)$ by creation operators.
The covariance condition requires that each $(\tau_i,\pi)$ is covariant, $i=1,2$. This means that $\psi_i(\varphi_i(a))=\pi(a)$ for $a\in \varphi_i^{-1}({\mathcal K}(E_i)), i=1,2$, where $\psi_i(\theta_{\xi,\eta})=\tau_i(\xi)\tau_i(\eta)^*$.
We get that 
\[{\mathcal T}_{E_2\otimes_A{\mathcal T}_{E_1}}\cong {\mathcal T}_{E_1\otimes_A{\mathcal T}_{E_2}}\cong {\mathcal T}_E\;\;{\mbox{and}} \;\;{\mathcal O}_{E_2\otimes_A{\mathcal O}_{E_1}}\cong {\mathcal O}_{E_1\otimes_A{\mathcal O}_{E_2}}\cong {\mathcal O}_E.\]

\medskip

{\bf 3.4 Example}. Let $A$ be a unital $C^*$-algebra, and $\alpha, \beta$ two commuting automorphisms of $A$. We assume that $\alpha$ and $\beta$  generate ${\mathbb Z}^2$ as a subgroup of $Aut(A)$. Denote by $A(\alpha), A(\beta)$ the Hilbert bimodules as in example 1.2.3, which are independent. We have $A(\alpha)\otimes_AA(\beta)\cong A(\beta)\otimes_AA(\alpha)$ by the map $\chi(\hat{a}\otimes\hat{b})=\widehat{\beta(\alpha^{-1}(a))}\otimes\hat{b}$. Indeed,
\[\langle\chi(\hat{a}_1\otimes\hat{b}_1), \chi(\hat{a}_2\otimes\hat{b}_2)\rangle=\langle\widehat{\beta(\alpha^{-1}(a_1))}\otimes\hat{b}_1,\widehat{\beta(\alpha^{-1}(a_2))}\otimes\hat{b}_2\rangle=\langle \hat{b}_1,\langle\widehat{\beta(\alpha^{-1}(a_1))},\widehat{\beta(\alpha^{-1}(a_2))}\rangle\cdot \hat{b}_2\rangle=\]
\[b_1^*\alpha(\beta(\alpha^{-1}(a_1))^*\beta(\alpha^{-1}(a_2)))b_2=b_1^*\beta(a_1^*a_2)b_2=\langle \hat{b}_1,\langle\hat{a}_1,\hat{a}_2\rangle\cdot\hat{b}_2\rangle=\langle\hat{a}_1\otimes\hat{b}_1,\hat{a}_2\otimes\hat{b}_2\rangle,\]
and
\[\chi(a'\cdot(\hat{a}\otimes\hat{b})\cdot b')=\chi(\widehat{\alpha(a')a}\otimes\widehat{bb'})=\widehat{\beta(\alpha^{-1}(\alpha(a')a))}\otimes\widehat{bb'}=\]\[\widehat{\beta(a'\alpha^{-1}(a))}\otimes\widehat{bb'}=\widehat{\beta(a')\beta(\alpha^{-1}(a))}\otimes\widehat{bb'}=a'\cdot\chi(\hat{a}\otimes\hat{b})\cdot b'.\]
Denote by ${\mathcal T}_{\alpha}, {\mathcal T}_{\beta}$ the corresponding Toeplitz algebras, each generated by $A$ and an isometry as in example 3, \S 1. Then $A(\beta)\otimes_A{\mathcal T}_{\alpha}$ becomes
a Hilbert bimodule over ${\mathcal T}_{\alpha}$, isomorphic to ${\mathcal T}_{\alpha}(\beta)$, where $\beta$ is extended to ${\mathcal T}_{\alpha}$ in the natural way, fixing the isometry. 
Similarly, $A(\alpha)\otimes_A{\mathcal T}_{\beta}\simeq {\mathcal T}_{\beta}(\alpha)$. We have
\[{\mathcal T}_{ {\mathcal T}_{\alpha}(\beta)}\cong {\mathcal T}_{{\mathcal T}_{\beta}(\alpha)}\cong {\mathcal T}_E,\]
where $E$ is the product system over ${\mathbb N}^2$ constructed from $A(\alpha)$, $A(\beta)$, and the isomorphism $\chi$.
Also, we may consider the Hilbert bimodules $A(\beta)\otimes_A(A\rtimes_{\alpha}{\mathbb Z})\cong (A\rtimes_{\alpha}{\mathbb Z})(\beta)$, and $A(\alpha)\otimes_A(A\rtimes_{\beta}{\mathbb Z})\cong (A\rtimes_{\beta}{\mathbb Z})(\alpha)$, where again $\alpha$ and $\beta$ are extended to $A\rtimes_{\beta}{\mathbb Z}$ and $A\rtimes_{\alpha}{\mathbb Z}$, respectively, by fixing the unitaries implementing the actions of ${\mathbb Z}$. We have
\[{\mathcal T}_{(A\rtimes_{\alpha}{\mathbb Z})(\beta)}\cong {\mathcal T}_{\beta}\rtimes_{\alpha}{\mathbb Z}, \;\; {\mathcal T}_{(A\rtimes_{\beta}{\mathbb Z})(\alpha)}\cong {\mathcal T}_{\alpha}\rtimes_{\beta}{\mathbb Z},\]
\[{\mathcal O}_{(A\rtimes_{\alpha}{\mathbb Z})(\beta)}\cong {\mathcal O}_{(A\rtimes_{\beta}{\mathbb Z})(\alpha)}\cong {\mathcal O}_E\cong A\rtimes_{\alpha,\beta}{\mathbb Z}^2.\]

\medskip

{\bf 3.5 Example}. Let $E=F={\mathbb C}^2$ with the usual structures of Hilbert bimodules over $A={\mathbb C}$. Then ${\mathcal O}_{E}\cong {\mathcal O}_2$, and $F\otimes {\mathcal O}_2$ becomes a Hilbert module over ${\mathcal O}_2$ with the usual operations. The left action will depend on a fixed isomorphism $\chi:E\otimes F\rightarrow F\otimes E$. If $\{e_1, e_2\}$ and $\{f_1, f_2\}$ are the canonical bases in $E$ and $F$ respectively, and $\chi_1(e_i\otimes f_j)=f_j\otimes e_i, \; i,j=1,2$, then the left multiplication is given by 
$S_i(f_j\otimes S_k)=f_j\otimes S_iS_k$, where ${\mathcal O}_2$ has generators $\{S_1,S_2\}$, and ${\mathcal O}_{F\otimes^{\chi_1} {\mathcal O}_2}\cong {\mathcal O}_2\otimes {\mathcal O}_2\cong{\mathcal O}_2$. On the other hand, if $\chi_2(e_i\otimes f_j)=f_i\otimes e_j,\; i,j=1,2$, then the left action is given by
$S_i(f_j\otimes S_k)=f_i\otimes S_jS_k$, and the ${\mathcal O}_2$-bimodule $F\otimes^{\chi_2}{\mathcal O}_2$ is degenerate. The corresponding Cuntz-Pimsner algebra ${\mathcal O}_{F\otimes^{\chi_2} {\mathcal O}_2}$ is isomorphic to $C({\mathbb T})\otimes {\mathcal O}_2$.
Indeed, if we interpret $E$ and $F$ as being each associated to the $1$-graph $\Gamma$ defining ${\mathcal O}_2$ ($\Gamma$ has two edges and one vertex), then the isomorphism
$\chi_2$ is defining the rank $2$ graph $\phi^*(\Gamma)$ where $\phi:{\mathbb N}^2\to{\mathbb N}, \phi(m,n)=m+n$ , and the last assertion follows from Example 6.1 and Proposition 2.10 in \cite{KP}. 

Note that  an arbitrary isomorphism $\chi$ does not necessarily define a rank 2 graph. For example, $\chi$ could be given by the unitary matrix $U=(u_{kl})$, where
\[u_{11}=\cos\alpha, u_{12}=-\sin\alpha, u_{21}=\sin\alpha, u_{22}=\cos\alpha,\]
\[u_{33}=\cos\beta, u_{34}=-\sin\beta, u_{43}=\sin\beta, u_{44}=\cos\beta,\]
for some angles $\alpha$ and $\beta$, and the rest of the entries equal to zero. For other  $C^*$-algebras defined by a product system of finite dimensional Hilbert spaces over the semigroup ${\mathbb N}^2$, see \cite{F1}.

\medskip

{\bf 3.6 Example}. Consider two star-commuting onto local homeomorphisms $\sigma_1$ and $\sigma_2$ of a compact space $X$. By definition, that means  that for every $x,y\in X$ such that $\sigma_1(x)=\sigma_2(y)$, there exists a unique $z\in X$ such that $\sigma_2(z)=x$ and $\sigma_1(z)=y$. This condition ensures that the associated conditional expectations $P_1, P_2$ commute (see \cite{ER}). Unfortunately, this condition was omitted in the Proposition on page 8 in \cite{D2}. I would like to thank Ruy and Jean for pointing this to me.

We can define the Hilbert bimodules $E_i=A(\alpha_i,P_i), i=1,2$ over $A=C(X)$ as in example 1.2.5.  Then there is an isomorphism $h:A(\alpha_1,P_1)\otimes_AA(\alpha_2,P_2)\to A(\alpha_1\circ\alpha_2,P_1\circ P_2)$, given by $h(\hat{a}\otimes\hat{b})=\widehat{a\alpha_1(b)}$, with inverse $h^{-1}(\hat{x})=\hat{x}\otimes\hat{1}$, which induces an isomorphism $\chi:E_1\otimes_AE_2\to E_2\otimes_AE_1$.
The resulting $C^*$-algebra ${\mathcal O}_{E_2\otimes_A{\mathcal O}_{E_1}}\cong {\mathcal O}_{E_1\otimes_A{\mathcal O}_{E_2}}$ can be understood as a crossed product of $C(X)$ by the semigroup ${\mathbb N}^2$. For other examples of semigroups of local homeomorphisms, see \cite{ER}.
\bigskip

\centerline{\S 4. {\small EXACT SEQUENCES IN $K$-THEORY}}

\bigskip

To study the $C^*$-algebra of Toeplitz operators on the quarter plane, Douglas and Howe (see \cite{DH}) considered the commutative diagram with exact rows and columns, where $j$ is the inclusion map, and $\pi$ is the quotient map:
\[\begin{array}{ccccccccc}{}&{}&0&{}&0&{}&0&{}&{}\\{}&{}&\downarrow &{}&\downarrow &{}&\downarrow &{}&{}\\0&\rightarrow &{\mathcal K}\otimes{\mathcal K}&\stackrel{j\otimes 1}{\longrightarrow} &{\mathcal T}\otimes{\mathcal K}&\stackrel{\pi\otimes 1}{\longrightarrow} &C({\mathbb T})\otimes{\mathcal K}&\rightarrow &0\\{}&{}&{}^{1\otimes j}\downarrow {\hspace {7mm}}&{}&{}^{1\otimes j}\downarrow\hspace{7mm} &{}&{}^{1\otimes j}\downarrow\hspace{7mm} &{}&{}\\0&\rightarrow &{\mathcal K}\otimes{\mathcal T}&\stackrel{j\otimes 1}{\longrightarrow} &{\mathcal T}\otimes{\mathcal T}&\stackrel{\pi\otimes 1}{\longrightarrow} &C({\mathbb T})\otimes{\mathcal T}&\rightarrow &0\\{}&{}&{}^{1\otimes\pi}\downarrow\hspace{7mm} &{}&{}^{1\otimes\pi}\downarrow\hspace{7mm} &{}&{}^{1\otimes\pi}\downarrow\hspace{7mm} &{}&{}\\0&\rightarrow &{\mathcal K}\otimes C({\mathbb T})&\stackrel{j\otimes 1}{\longrightarrow }&{\mathcal T}\otimes C({\mathbb T})&\stackrel{\pi\otimes 1}{\longrightarrow} &C({\mathbb T})\otimes C({\mathbb T})&\rightarrow &0\\{}&{}&\downarrow &{}&\downarrow &{}&\downarrow &{}&{}\\{}&{}&0&{}&0&{}&0&{}&{}\end{array}.\]

\medskip

{\bf 4.1 Corollary}.  We have the short exact sequences
\[0\rightarrow{\mathcal T}\otimes{\mathcal K}+{\mathcal K}\otimes{\mathcal T}\stackrel{1\otimes j+j\otimes 1}{\longrightarrow}{\mathcal T}\otimes{\mathcal T}\stackrel{\pi\otimes\pi}{\longrightarrow}C({\mathbb T})\otimes C({\mathbb T})\rightarrow 0,\]
\[0\rightarrow{\mathcal K}\otimes{\mathcal K}\stackrel{j\otimes 1+1\otimes j}{\longrightarrow}{\mathcal T}\otimes{\mathcal K}+{\mathcal K}\otimes{\mathcal T}\stackrel{\pi\otimes 1+1\otimes \pi}{\longrightarrow}C({\mathbb T})\otimes{\mathcal K}\oplus{\mathcal K}\otimes C({\mathbb T})\rightarrow 0.\]

The above $3\times 3$ diagram and the exact sequences in the corollary are particular cases of a more general situation, for which we provide a proof. 
 
\medskip
 
{\bf 4.2 Lemma}. Let $A$ be a $C^*$-algebra and $I, J$ two closed two-sided ideals of $A$. Then we have the commutative diagram with exact rows and columns, where the maps are the canonical ones:
\[\begin{array}{ccccccccc}{}&{}&0&{}&0&{}&0&{}&{}\\{}&{}&\downarrow &{}&\downarrow &{}&\downarrow &{}&{}\\0&\rightarrow &I\cap J&\stackrel{\lambda_1}{\longrightarrow} &I&\stackrel{\omega_1}{\longrightarrow} &I/(I\cap J)&\rightarrow &0\\{}&{}&{}^{\lambda_2}\downarrow\hspace{3mm} &{}&{}^{\iota_1}\downarrow\hspace{3mm} &{}&{}^{\sigma_1}\downarrow\hspace{3mm} &{}&{}\\0&\rightarrow  &J&\stackrel{\iota_2}{\longrightarrow} &A&\stackrel{\pi_2}{\longrightarrow}&A/J&\rightarrow&0\\{}&{}&{}^{\omega_2}\downarrow\hspace{3mm} &{}&{}^{\pi_1}\downarrow\hspace{3mm} &{}&{}^{\rho_2}\downarrow\hspace{3mm} &{}&{}\\0&\rightarrow &J/(I\cap J)&\stackrel{\sigma_2}{\longrightarrow} &A/I&\stackrel{\rho_1}{\longrightarrow} &A/(I+J)&\rightarrow &0\\{}&{}&\downarrow &{}&\downarrow &{}&\downarrow &{}&{}\\{}&{}&0&{}&0&{}&0&{}&{}\end{array}.\]
From this we get the exact sequences
\[0\rightarrow I+J\stackrel{\iota_1+\iota_2}{\longrightarrow} A\stackrel{\pi}{\longrightarrow} A/(I+J)\rightarrow 0,\]where $\pi=\rho_1\circ\pi_1=\rho_2\circ\pi_2$, and\[0\rightarrow I\cap J\stackrel{\lambda_1+\lambda_2}{\longrightarrow} I+J\stackrel{\omega_1+\omega_2}{\longrightarrow} I/(I\cap J)\oplus J/(I\cap J)\rightarrow 0.\]
Applying the $K$-theory functor, we get
\[K_0(I+J)\rightarrow K_0(A)\rightarrow K_0(A/(I+J))\]
\[\uparrow\hspace{40mm}\downarrow\]\[K_1(A/(I+J))\leftarrow K_1(A)\leftarrow K_1(I+J),\]
\bigskip
\[K_0(I\cap J)\rightarrow K_0(I+J)\rightarrow K_0(I/(I\cap J))\oplus K_0(J/(I\cap J))\]\[\uparrow\hspace{68mm}\downarrow\]\[K_1(I/(I\cap J))\oplus K_1(J/(I\cap J)\leftarrow K_1(I+J)\leftarrow K_1(I\cap J).\]

{\em Proof}. Note that the first two rows and the first two columns are obviously exact. For the third row, the map $\pi_1\circ\iota_2$ has kernel $I\cap J$. This defines a map $\sigma_2:J/(I\cap J)\to A/I$ such that $\sigma_2\circ\omega_2=\pi_1\circ\iota_2$. By the second isomorphism theorem, $(I+J)/I\cong J/(I\cap J)$, and by the third isomorphism theorem, $(A/I)/((I+J)/I)\cong A/(I+J)$, hence the third row is exact. The exactness of the third column is proved similarly. Consider now the diagonal morphism $\pi=\rho_1\circ\pi_1=\rho_2\circ\pi_2: A\rightarrow A/(I+J)$. If $a\in \ker\pi$, then $\pi_2(a)\in \ker\rho_2=\sigma_1(I/(I\cap J))=\sigma_1(\omega_1(I))$, hence there is $b\in I$ with $\pi_2(a)=\sigma_1(\omega_1(b))=\pi_2(\iota_1(b))$. It follows that $a-\iota_1(b)\in\ker\pi_2=\iota_2(J)$, and there is $c\in J$ with $a=\iota_1(b)+\iota_2(c)$. This gives the first exact sequence. For the second, we use the map $\omega_1+\omega_2:I+J\to I/(I\cap J)\oplus J/(I\cap J)$ which has kernel $I\cap J$. $\Box$

We  generalize the above diagram of Douglas and Howe to certain iterated Toeplitz and Cuntz-Pimsner algebras. By the lemma, we  get some exact sequences which we hope will help  to do $K$-theory computations in some particular cases.

 \medskip

{\bf 4.3 Theorem}. Consider a $C^*$-algebra $A$ and two finitely generated Hilbert bimodules $E_1, E_2$ with a fixed isomorphism $\chi:E_1\otimes_A E_2\rightarrow E_2\otimes_A E_1$. We assume that $E_2\otimes_A {\mathcal T}_{E_1}, E_2\otimes_A{\mathcal O}_{E_1}, E_1\otimes_A {\mathcal T}_{E_2}, E_1\otimes_A{\mathcal O}_{E_2}$ are nondegenerate as Hilbert bimodules, with the structure described in \S 3. Then we have the following commuting diagram (depending on $\chi$), with exact rows and columns, where the maps are canonical: 
\[\begin{array}{ccccccccc}{}&{}&0&{}&0&{}&0&{}&{}\\{}&{}&\downarrow &{}&\downarrow &{}&\downarrow &{}&{}\\0&\rightarrow &{\mathcal K}(\ell^2(E_1\otimes_A E_2))&\rightarrow &{\mathcal K}(\ell^2(E_2\otimes_A{\mathcal T}_{E_1}))&\rightarrow &{\mathcal K}(\ell^2(E_2\otimes_A{\mathcal O}_{E_1}))&\rightarrow &0\\{}&{}&\downarrow &{}&\downarrow &{}&\downarrow &{}&{}\\0&\rightarrow &{\mathcal K}(\ell^2(E_1\otimes_A{\mathcal T}_{E_2}))&\rightarrow &{\mathcal T}_{E_1\otimes_A{\mathcal T}_{E_2}}\cong {\mathcal T}_{E_2\otimes_A{\mathcal T}_{E_1}}&\rightarrow &{\mathcal O}_{E_1\otimes_A{\mathcal T}_{E_2}}\cong{\mathcal T}_{E_2\otimes_A{\mathcal O}_{E_1}}&\rightarrow &0\\{}&{}&\downarrow &{}&\downarrow &{}&\downarrow &{}&{}\\0&\rightarrow &{\mathcal K}(\ell^2(E_1\otimes_A{\mathcal O}_{E_2}))&\rightarrow &{\mathcal T}_{E_1\otimes_A{\mathcal O}_{E_2}}\cong {\mathcal O}_{E_2\otimes_A{\mathcal T}_{E_1}}&\rightarrow &{\mathcal O}_{E_1\otimes_A{\mathcal O}_{E_2}}\cong{\mathcal O}_{E_2\otimes_A{\mathcal O}_{E_1}}&\rightarrow &0\\
{}&{}&\downarrow &{}&\downarrow &{}&\downarrow &{}&{}\\{}&{}&0&{}&0&{}&0&{}&{}\end{array}.\]

{\em Proof}. Recall that the algebra ${\mathcal T}_{E_i}$ is represented on the Fock space $\ell^2(E_i)$ and we have a short exact sequence
\[0\rightarrow {\mathcal K}(\ell^2(E_i))\rightarrow {\mathcal T}_{E_i}\rightarrow {\mathcal O}_{E_i}\rightarrow 0,\]
for $i=1,2$. We apply the above lemma for the $C^*$-algebra ${\mathcal T}_{E_1\otimes_A{\mathcal T}_{E_2}}\cong {\mathcal T}_{E_2\otimes_A{\mathcal T}_{E_1}}$ with ideals $I={\mathcal K}(\ell^2(E_2\otimes_A{\mathcal T}_{E_1}))$ and $J= {\mathcal K}(\ell^2(E_1\otimes_A{\mathcal T}_{E_2}))$. The nondegeneracy assumption implies that $\ell^2(E_2\otimes_A{\mathcal T}_{E_1})\cong \ell^2(E_2)\otimes_A{\mathcal T}_{E_1}$ and $\ell^2(E_2\otimes_A{\mathcal O}_{E_1})\cong \ell^2(E_2)\otimes_A{\mathcal O}_{E_1}$. Note that the map $\chi:E_1\otimes_A E_2\to E_2\otimes_AE_1$ induces  isomorphisms ${\mathcal K}(\ell^2(E_2\otimes_A{\mathcal K}(\ell^2(E_1))))\cong {\mathcal K}(\ell^2(E_1\otimes_A{\mathcal K}(\ell^2(E_2))))\cong {\mathcal K}(\ell^2(E_1\otimes E_2))$, therefore $I\cap J={\mathcal K}(\ell^2(E_1\otimes_A E_2))$. $\Box$

\medskip

{\bf 4.4 Corollary}. Under the same assumptions as in the theorem,  we get the short exact sequences
\[0\rightarrow{\mathcal K}(\ell^2(E_2\otimes_A{\mathcal T}_{E_1}))+{\mathcal K}(\ell^2(E_1\otimes_A{\mathcal T}_{E_2}))\rightarrow {\mathcal T}_{E_1\otimes_A{\mathcal T}_{E_2}}\rightarrow{\mathcal O}_{E_1\otimes_A{\mathcal O}_{E_2}}\rightarrow 0\]
and
\[0\rightarrow{\mathcal K}(\ell^2(E_1\otimes_AE_2))\rightarrow{\mathcal K}(\ell^2(E_2\otimes_A{\mathcal T}_{E_1}))+{\mathcal K}(\ell^2(E_1\otimes_A{\mathcal T}_{E_2}))\rightarrow\]\[\rightarrow{\mathcal K}(\ell^2(E_2\otimes_A{\mathcal O}_{E_1}))\oplus{\mathcal K}(\ell^2(E_1\otimes_A{\mathcal O}_{E_2}))\rightarrow 0.\]
The corresponding six-term exact sequences of $K$-theory give us information about the $K$-theory of ${\mathcal O}_{E_1\otimes_A{\mathcal O}_{E_2}}$, once we identify the maps between the various $K$-groups.
Note that the  Toeplitz algebras are $KK$-equivalent to the  $C^*$-algebra $A$.

\medskip

{\bf 4.5 Example}. For $\alpha_1,\alpha_2$ two commuting independent automorphisms of a $C^*$-algebra $A$, the diagram is
\[\begin{array}{ccccccccc}{}&{}&0&{}&0&{}&0&{}&{}\\{}&{}&\downarrow &{}&\downarrow &{}&\downarrow &{}&{}\\0&\rightarrow &{\mathcal K}\otimes  A\otimes{\mathcal K}&\rightarrow &{\mathcal T}_{\alpha_1}\otimes{\mathcal K}&\rightarrow &A\rtimes_{\alpha_1}{\mathbb Z}\otimes{\mathcal K}&\rightarrow &0\\{}&{}&\downarrow &{}&\downarrow &{}&\downarrow &{}&{}\\0&\rightarrow &{\mathcal K}\otimes{\mathcal T}_{\alpha_2}&\rightarrow &{\mathcal T}_{{\mathcal T}_{\alpha_2}(\alpha_1)}\cong {\mathcal T}_{{\mathcal T}_{\alpha_1}(\alpha_2)}&\rightarrow &{\mathcal T}_{\alpha_2}\rtimes_{\alpha_1}{\mathbb Z}\cong{\mathcal T}_{(A\rtimes_{\alpha_1}{\mathbb Z})(\alpha_2)}&\rightarrow &0\\{}&{}&\downarrow &{}&\downarrow &{}&\downarrow &{}&{}\\0&\rightarrow &{\mathcal K}\otimes A\rtimes_{\alpha_2}{\mathbb Z}&\rightarrow &{\mathcal T}_{\alpha_1}\rtimes_{\alpha_2}{\mathbb Z}\cong{\mathcal T}_{(A\rtimes_{\alpha_2}{\mathbb Z})(\alpha_1) }&\rightarrow &A\rtimes_{\alpha_1,\alpha_2}{\mathbb Z}^2&\rightarrow &0\\
{}&{}&\downarrow &{}&\downarrow &{}&\downarrow &{}&{}\\{}&{}&0&{}&0&{}&0&{}&{}\end{array}.\]

We get the exact sequences
\[K_0(A)\rightarrow K_0({\mathcal T}_{\alpha_1}\otimes{\mathcal K}+{\mathcal K}\otimes{\mathcal T}_{\alpha_2})\rightarrow K_0(A\rtimes_{\alpha_1}{\mathbb Z})\oplus K_0(A\rtimes_{\alpha_2}{\mathbb Z})\]\[\uparrow\partial^1+\partial^2\hspace{65mm}\partial^1+\partial^2\downarrow\]
\[K_1(A\rtimes_{\alpha_1}{\mathbb Z})\oplus K_1(A\rtimes_{\alpha_2}{\mathbb Z})\leftarrow K_1({\mathcal T}_{\alpha_1}\otimes{\mathcal K}+{\mathcal K}\otimes{\mathcal T}_{\alpha_2})\leftarrow K_1(A)\]
\[K_0({\mathcal T}_{\alpha_1}\otimes{\mathcal K}+{\mathcal K}\otimes{\mathcal T}_{\alpha_2})\rightarrow K_0(A)\rightarrow K_0(A\rtimes_{\alpha_1,\alpha_2}{\mathbb Z}^2)\]
\[\uparrow\hspace{55mm}\downarrow\]
\[K_1(A\rtimes_{\alpha_1,\alpha_2}{\mathbb Z}^2)\leftarrow K_1(A)\leftarrow K_1({\mathcal T}_{\alpha_1}\otimes{\mathcal K}+{\mathcal K}\otimes{\mathcal T}_{\alpha_2}).\]

For $A=C(X)$ with $X$ a Cantor set, and $\sigma,\tau:X\rightarrow X$ two commuting  local homeomorphisms. Denote by  $\alpha, \beta$  the induced endomorphisms of $A$, and assume that the associated conditional expectations commute. Applying the theorem for the Hilbert bimodules $E=A(\alpha, P)$, $F=A(\beta, Q)$ and the canonical isomorphism $E\otimes_AF\cong F\otimes_AE$, we get
\[C(X,{\mathbb Z})\rightarrow K_0({\mathcal T}_{E}\otimes {\mathcal K}+{\mathcal K}\otimes{\mathcal T}_{F})\rightarrow C(X,{\mathbb Z})/im(1-\alpha_*)\oplus C(X,{\mathbb Z})/im(1-\beta_*)\]
\[\uparrow\hspace{74mm}\downarrow\hspace{37mm}\]
\[ker(1-\alpha_*)\oplus ker(1-\beta_*)\leftarrow K_1({\mathcal T}_{E}\otimes {\mathcal K}+{\mathcal K}\otimes{\mathcal T}_{F})\leftarrow 0\hspace{46mm}\]

\[K_0({\mathcal T}_{E}\otimes {\mathcal K}+{\mathcal K}\otimes{\mathcal T}_{F})\rightarrow C(X,{\mathbb Z})\rightarrow K_0({\mathcal O}_{F\otimes_A{{\mathcal O}_E}})\]
\[\uparrow\hspace{50mm}\downarrow\]
\[K_1({\mathcal O}_{F\otimes_A{{\mathcal O}_E}})\leftarrow 0\leftarrow K_1({\mathcal T}_{E}\otimes {\mathcal K}+{\mathcal K}\otimes{\mathcal T}_{F}).\]
Here we used the fact that $K_0(A)=C(X,{\mathbb Z})$ and $K_1(A)=0$.

\medskip

{\bf 4.6 Example}. Let $E={\mathbb C}^m, F={\mathbb C}^n$ and fix $\chi:E\otimes F\rightarrow F\otimes E$ an isomorphism. The corresponding diagram (depending on $\chi$) is
\[\begin{array}{ccccccccc}{}&{}&0&{}&0&{}&0&{}&{}\\{}&{}&\downarrow &{}&\downarrow &{}&\downarrow &{}&{}\\0&\rightarrow &{\mathcal K}\otimes{\mathcal K}&\rightarrow &{\mathcal E}_m\otimes{\mathcal K}&\rightarrow &{\mathcal O}_m\otimes{\mathcal K}&\rightarrow &0\\{}&{}&\downarrow &{}&\downarrow &{}&\downarrow &{}&{}\\0&\rightarrow &{\mathcal K}\otimes{\mathcal E}_n&\rightarrow &{\mathcal T}_{{\mathbb C}^m\otimes{\mathcal E}_n}\cong {\mathcal T}_{{\mathbb C}^n\otimes{\mathcal E}_m}&\rightarrow &{\mathcal O}_{{\mathbb C}^m\otimes{\mathcal E}_n}\cong{\mathcal T}_{{\mathbb C}^n\otimes{\mathcal O}_m}&\rightarrow &0\\{}&{}&\downarrow &{}&\downarrow &{}&\downarrow &{}&{}\\0&\rightarrow &{\mathcal K}\otimes {\mathcal O}_n&\rightarrow &{\mathcal T}_{{\mathbb C}^m\otimes{\mathcal O}_n}\cong{\mathcal O}_{{\mathbb C}^n\otimes{\mathcal E}_m} &\rightarrow &{\mathcal O}_{{\mathbb C}^m\otimes{\mathcal O}_n}\cong{\mathcal O}_{{\mathbb C}^n\otimes{\mathcal O}_m}&\rightarrow &0\\{}&{}&\downarrow &{}&\downarrow &{}&\downarrow &{}&{}\\{}&{}&0&{}&0&{}&0&{}&{}\end{array}.\]

The exact sequences (depending on $\chi$) are
\[{\mathbb Z}\rightarrow K_0({\mathcal E}_m\otimes{\mathcal K}+{\mathcal K}\otimes{\mathcal E}_n)\rightarrow{\mathbb Z}_{m-1}\oplus{\mathbb Z}_{n-1}\]\[\uparrow\hspace{53mm}\downarrow\hspace{10mm}\]\[0\;\longleftarrow \;K_1({\mathcal E}_m\otimes{\mathcal K}+{\mathcal K}\otimes{\mathcal E}_n)\;\longleftarrow\; 0\hspace{10mm}\]
and
\[K_0({\mathcal E}_m\otimes{\mathcal K}+{\mathcal K}\otimes{\mathcal E}_n)\rightarrow{\mathbb Z}\rightarrow K_0({\mathcal O}_{{\mathbb C}^m\otimes{\mathcal O}_n})\]
\[\hspace{18mm}\uparrow\hspace{36mm}\downarrow\hspace{10mm}\]
\[K_1({\mathcal O}_{{\mathbb C}^m\otimes{\mathcal O}_n})\hspace{2mm}\longleftarrow \hspace{2mm}0\hspace{2mm}\longleftarrow\hspace{3mm}
0.\]

In particular, if $\chi$ is just the flip, we have
\[{\mathcal T}_{{\mathbb C}^m\otimes{\mathcal E}_n}\cong {\mathcal T}_{{\mathbb C}^n\otimes{\mathcal E}_m}\cong {\mathcal E}_n\otimes{\mathcal E}_m, \;\; {\mathcal O}_{{\mathbb C}^m\otimes{\mathcal E}_n}\cong{\mathcal T}_{{\mathbb C}^n\otimes{\mathcal O}_m}\cong {\mathcal O}_m\otimes{\mathcal E}_n,\]
\[{\mathcal T}_{{\mathbb C}^m\otimes{\mathcal O}_n}\cong{\mathcal O}_{{\mathbb C}^n\otimes{\mathcal E}_m}\cong {\mathcal E}_m\otimes{\mathcal O}_n, \;\; {\mathcal O}_{{\mathbb C}^m\otimes{\mathcal O}_n}\cong{\mathcal O}_{{\mathbb C}^n\otimes{\mathcal O}_m}\cong {\mathcal O}_m\otimes{\mathcal O}_n,\]
and we recover the $K$-theory of $ {\mathcal O}_m\otimes{\mathcal O}_n$.

\medskip

{\bf 4.7 Example}. Let $A=C({\mathbb T})$ and $\sigma_i(x)=x^{p_i}, i=1,2$ with $p_1, p_2$ relatively prime. Then the associated conditional expectations commute, and using
the $K$-theory computations in \cite{D1}, we get
\[{\mathbb Z}\rightarrow K_0({\mathcal T}_{\alpha_1}\otimes{\mathcal K}+{\mathcal K}\otimes{\mathcal T}_{\alpha_2})\rightarrow{\mathbb Z}^2\oplus{\mathbb Z}_{p_1-1}\oplus{\mathbb Z}_{p_2-1}\]
\[\uparrow\hspace{73mm}\downarrow\hspace{3mm}\]
\[ {\mathbb Z}\;\oplus \;{\mathbb Z}\hspace{4mm}\longleftarrow\hspace{4mm} K_1({\mathcal T}_{\alpha_1}\otimes{\mathcal K}+{\mathcal K}\otimes{\mathcal T}_{\alpha_2})\hspace{4mm}\longleftarrow\hspace{4mm} {\mathbb Z}\hspace{7mm}\]

\[K_0({\mathcal T}_{\alpha_1}\otimes{\mathcal K}+{\mathcal K}\otimes{\mathcal T}_{\alpha_2})\rightarrow {\mathbb Z}\rightarrow K_0({\mathcal O}_{E_2\otimes_A{\mathcal O}_{E_1}})\]
\[\uparrow\hspace{50mm}\downarrow\]
\[K_1({\mathcal O}_{E_2\otimes_A{\mathcal O}_{E_1}})\leftarrow {\mathbb Z}\leftarrow K_1({\mathcal T}_{\alpha_1}\otimes{\mathcal K}+{\mathcal K}\otimes{\mathcal T}_{\alpha_2}).\]
More generally, we may consider coverings of the n-torus ${\mathbb T}^n$.


\end{document}